\documentclass[11pt, leqno]{amsart}
\usepackage{amsmath}
\usepackage{amsfonts,amssymb,amsthm,amstext,amscd,boxedminipage}  %
\usepackage[mathscr]{eucal}

\usepackage{graphicx,psfrag} 
\usepackage{color}
 \numberwithin{equation}{section}

\renewcommand{\setminus}{{\smallsetminus}}  %
\DeclareMathSymbol{\minus} {\mathord}{operators}{"2D} %

\newtheorem{theo}{Theorem}[section]

\newtheorem{lem}[theo]{Lemma}
\newtheorem{defi}[theo]{Definition}

\newtheorem{corol}[theo]{Corollary}

\def \proof {{\bf Proof$\colon$}\ }
\def \Z{{\mathbb {Z}}}

\def \N{{\mathbb {N}}}

\def\nt{$n$-trivializer}

\def\no{\noindent}
\begin{document}

\title[Quantum invariants and volume]{A note on quantum 3-manifold invariants and hyperbolic volume}

\author[E. Kalfagianni]{Efstratia Kalfagianni}

\address[]{Michigan State
University,
E. Lansing, MI, 48823}
\email[]{kalfagia@@math.msu.edu}

\begin{abstract} {For a closed, oriented 3-manifold $M$ and an integer $r>0$,  let $\tau_r(M)$ denote 
the $SU(2)$
Reshetikhin-Turaev-Witten invariant of $M$, at level $r$. We show that for every 
$n>0$, and for $r_1, \ldots, r_n >0$ sufficiently large integers, there exist infinitely many non-homeomorphic
hyperbolic 
3-manifolds
$M$, all of which have different hyperbolic volume, and such that
$\tau_{r_i}(M)=1$,
for
$i=1,\ldots, n$.

\medskip
\smallskip
\medskip

\smallskip

\medskip

}
\end{abstract}

\maketitle
\smallskip

\smallskip

\section{introduction}  Throughout the paper, $M$ will denote a closed orientable 3-manifold.
The Reshetikhin-Turaev-Witten $SU(2)$-invariants  of $M$ are
complex valued numbers (\cite{rt}, \cite{kn:kim}) parametrized by positive integers (levels).
Given  a positive integer $r\in \N$, let $\tau_r(M)$ denote the
$SU(2)$ invariant of
$M$ at level $r$. In the case that $M$ is hyperbolic,
let 
${\rm vol }(M)$ denote the hyperbolic volume of $M$. It has been speculated
that ${\rm vol }(M)$ is determined by the entire collection of
the invariants $T_M:=\{\tau_r(M)| \ r=1,2,..\}$ (see \cite{kn:mu}). However, at the moment, it is not clear
what the precise statement of a conjecture in this direction should be.
In this paper we are concerned with the question of the extent to which the
finite sub sequences of $T_M$ determine ${\rm vol }(M)$.
The main result is the following theorem that shows that for most finite sub sequences
the answer to this question is an emphatic {\sl no}.

\begin{theo} \label{theo:main} Fix $n>0$. There is a constant $C_n>0$ such that
for every $n$-tuple of integers  $r_1, \ldots, r_n >C_n$, there exist 
an infinite sequence 
of hyperbolic 3-manifolds
$\{M_k\}_{k\in \N}$ such that 

\no (a) \ \ For
$i=1,\ldots, n$, and every $k\in \N$ we have  $\tau_{r_i}(M_k)=1$; and  

\no (b)\ \  ${\displaystyle{\ldots > {\rm vol }(M_{k+1})>{\rm vol }(M_k)>
\ldots {\rm vol }(M_0)>{n\over 2}}}$.
\end{theo}

Given a framed link $L$ in $S^3$, there is a sequence of Laurent polynomials $\{ J_N(L,t)\}_{N\in \N}$; the colored Jones
polynomials (\cite{kn:kim}). 
For the trivial knot $U$, if equipped with the 0-framing,  we have
$$J_N(U,t)=[N]:={{t^N-t^{-N}}\over {t-t^{-1}}}.$$ To prove
Theorem \ref{theo:main} we need the following:

\begin{theo} \label{theo:mainknot}Given integers $n>0$ and $r_1, \ldots, r_n >2$, 
there is a knot
$K\subset S^3$ such that for any common framing on $K$ and $U$ we have 

$$J_N(K, e_{r_i})=J_N(U, e_{r_i}) \ \ {\rm for \ all}\ \ N\in \N.$$

\no Here, for $i=1, \ldots, n$, $e_{r_i}:= e^{{2\pi {\sqrt {-1}}\over {r_i}}}$ is
a primitive $r_i$-th root of unity.
For fixed $n$, if $r_1, \ldots, r_n$ are sufficiently large, then
$K$ can be chosen hyperbolic with  ${\rm vol }(S^3\setminus K)>n$. Furthermore, if $M$ is a hyperbolic 3-manifold
obtained by ${\displaystyle {p\over q}}$-surgery on $K$, for some $|q|>12$, then we have
$${\rm vol}(M) \: \geq \: {(1-{{127}\over {q^2}})}^{3\over 2}\ n . \eqno(1)$$
\end{theo}
\smallskip

The proof of the first part of Theorem \ref{theo:main} uses a result of Lackenby (\cite{kn:la}) 
and a
construction of \cite{kn:k}. For the remaining claims, we need
Thurston's hyperbolic Dehn surgery theorem (\cite{kn:thurston}) and a result proved jointly with  Futer and  Purcell
(\cite{dej}).

\begin{corol} \label{corol:nonuniform} Given an integer $n>0$, there is a sequence of  hyperbolic 3-manifolds
$\{M_i\}_{i\in \N}$ and an increasing sequence of positive integers  $\{m_i\}_{i\in \N}$ such that
$${\rm vol}(M_i)>{n \over 2}\ \ {\rm  and } \ \ \tau_{m_i}(M_i)=1,$$ for all $i\in \N$.      
\end{corol}

\smallskip

\section{The Proofs}
\subsection{Some properties of the colored Jones polynomials}A crossing disc of a knot $J$
is an  embedded disc $D\subset S^3$
that intersects $J$ only in its interior 
exactly twice geometrically and with zero algebraic intersection number.
The curve $\partial D$ is a crossing circle
for $J$. A knot $K$ is said to be obtained from $J$ by a generalized crossing of order
$r \in \Z$ iff $K$  is the result of $J$ under surgery of $S^3$ along $\partial D$
with surgery slope ${\displaystyle {{ 1 \over r} }}$.

\begin{defi} { \rm ( Definition 1.1, \cite{kn:la}) }  \label{defi:fox} Let $r\in \N$ and let $J$ and $K$ be two
0-framed knots in $S^3$. We say $K$ and $J$ are 
congruent modulo $(r,\ 2)$, iff $K$ is obtained from $J$ by a collection of generalized
crossing changes of order $r$ supported on disjoint crossing discs.
In this case we will write $J\equiv K({\rm mod}(r, \ 2))$.
\end{defi}

For $j\in N$,  let $K^j$ denote the $j$-th parallel cable of $K$ formed with 0-
framing and let $J(K^j, t)$ denote the Jones polynomial of $K^j$. We need the
following result of Lackenby.

\begin{lem} \label{lem:lac} {\rm (Corollary 2.8, \cite{kn:la})} Let $r>2$ be an integer and  let $e_r:= e^{{2\pi {\sqrt {-1}}\over r}}$
denote a primitive $r$-th root of unity. Suppose that $J$ and $K$ are
0-framed knots in $S^3$.
If $J\equiv K({\rm mod}(r, \ 2))$, then,

$$J(K^j, e_r)=J({J}^j, e_r), \ \ {\rm for \ all}\ \ j\in \N.$$
\end{lem}

We recall that  for a 0-framed link $L$ the value $J_N(L, e_r)$ is a linear combination of $J(L^j, e_r)$,
with the coefficients of the combination being constants independent of $L$ (Theorem 4.15, \cite{kn:kim}).
Using this fact and  Lemma \ref{lem:lac} we have:

\begin{corol} \label{corol:colored} Let $r>2$ be an integer and
let $e_r:= e^{{2\pi {\sqrt {-1}}\over r}}$ denote
a primitive $r$-th root of unity. Suppose that $J$ and $K$ are
0-framed knots in $S^3$.
If $J\equiv K({\rm mod}(r, \ 2))$,
then, for every integer $q>0$, we have $ J_N(K^q, e_r)=J_N(J^q, e_r)$, for  all
$N\in \N.$
In particular, we have $  J_N(K, e_r)=J_N(J, e_r)$, 
for  all $N\in \N.$
\end{corol}
\smallskip

\subsection{ A construction of hyperbolic Brunnian links}
To prove Theorem \ref{theo:mainknot} we will need the following lemma which summarizes results proved in \cite{kn:k}
and uses results proved jointly with Askitas in \cite{kn:ak}.
Below we will  sketch the proof referring the reader to the original references for details.
\smallskip

\begin{lem} \label{lem:brunnian}
For every $n>0$, there is an $(n+1)$-component
link 
$L_n:= {U} \cup K_1 \cup\ldots K_n$ with the following properties:
\vskip 0.04in

\no (a) \  \ $L_n$ is Brunnian; that is every proper sublink of $L_n$ is a trivial link.
\vskip 0.04in

\no (b)\ \  For $i=1\ldots, n$, $K_i$ bounds a crossing disc $D_i\subset S^3$ of $U$.
\vskip 0.04in

\no (c)\ \ $L_n$ is hyperbolic; that is the interior of the 3-manifold ${\overline{M_n:=S^3\setminus \eta(L_n)}}$ admits a complete hyperbolic
metric of finite volume. Here, $\eta(L_n)$ denotes a tubular neighborhood of $L_n$.
\vskip 0.04in

\no (d)\ \ For $n>1$,  any  collection of generalized crossing changes along any collection of discs
formed by a proper subset of
$\{ D_1,\ldots, D_n\}$, leaves $U$ unknotted.
\vskip 0.04in

\no (e)\ \  Every knot obtained by a collection of $n$ generalized crossing changes of order
$r_1,\ldots, r_n>0$ along
$D_1, \ldots D_n$, respectively, is non-trivial.
\end{lem}
\proof For $n=1$, we take $L_1:=U\cup K_1$ to be the Whitehead link
and for $n=2$, we take $L_2:=U\cup K_1\cup K_2$ to be the 3-component Borromean link
(see Figure 1). It is well known that they are both hyperbolic (\cite{kn:thurston}).
For $n=2$, condition  (d) is clearly satisfied and for $n=1$ and $n=2$, condition (e) 
is true since the resulting knot will be a non-trivial twist knot.
\medskip

\begin{figure}[htbp] %
   \centering
   \includegraphics[width=3in, ]{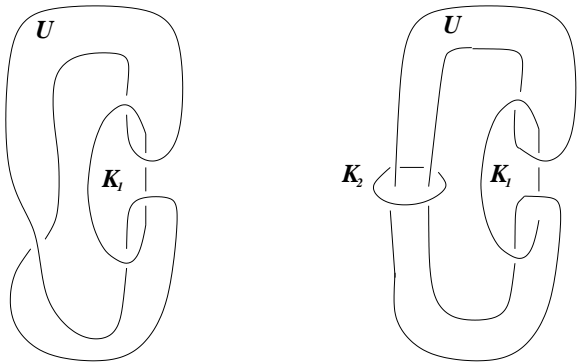} 
   
   \caption{The link $L_n$ for $n=1$ and $n=2$.}
  \end{figure}

For $n>2$ the construction of a link $L_n$ as claimed above, is given in 
the proof of Theorem 3.1 of  \cite{kn:k}. We now outline the construction using the terminology of
\cite{kn:ak} and  \cite{kn:k}. Recall that a knot ${\hat K}$ is $n$-adjacent to the unknot, for some
$n\in {\bf N}$,
if ${\hat K}$ admits an embedding containing $n$ generalized crossings
such that changing any $0<m\leq n$ of them yields an
embedding of  the unknot. A collection of crossing circles
corresponding to these crossings is called an \nt.
Theorem 3.5 of \cite{kn:k} proves the following:
For 
every $n>2$, there exists a knot ${\hat K}$ that is $n$-adjacent to the unknot, and it
admits an \nt \ $K_1 \cup\ldots\cup K_n$ such that the link $L^{*}_n={\hat K}\cup K_1 \cup\ldots \cup K_n$
is hyperbolic.
Let $D_1, \ldots, D_n$ denote crossing discs corresponding to $K_1, \ldots, K_n$, respectively.
Now let $L_{n}$ denote the $(n+1)$-component link obtained from $L^{*}_n$ by performing the $n$ generalized crossing
changes (along $D_1, \ldots, D_n$) that exhibit ${\hat K}$ as n-adjacent to the unknot, simultaneously.  
We can write $L_n:= {U} \cup K_1 \cup\ldots \cup K_n$, where $U$ is the unknot resulting from  ${\hat K}$
after these crossing changes. Since the links $L^{*}_n$ and $L_n$ differ by twists along a collection of discs,
they have 
homeomorphic complements. We conclude that $L_n$ is hyperbolic.
Thus $L_n$ satisfies conditions (b) and (c) of the statement of the lemma. 
Next, let us focus on a proper sub collection
of crossing discs; say, without loss of generality, $D_1, \ldots, D_{n-1}$. Since $K_1 \cup\ldots K_n$
is an \nt \ for ${\hat K}$, by the definition of adjacency to the unknot,   $K_1 \cup\ldots \cup K_{n-1}$ is an $(n-1)$-trivializer \ 
By Theorem 2.2 of \cite{kn:ak}, and its proof, $U$ bounds an embedded disc, say $\Delta$, in the complement of $K_1 \cup\ldots\cup K_{n-1}$
such that $D_i\cap \Delta$ is a single arc properly embedded on $\Delta$. It follows that $U\cup K_1 \cup\ldots\cup  K_{n-1}$
is the trivial link and that every collection of generalized crossing changes supported along $D_1, \ldots, D_{n-1}$
leaves $U$ unknotted. This proves (b) and (d). To see  (e) suppose, on the contrary, that there is a collection of $n$ generalized crossing changes of order
$r_1,\ldots, r_n>0$ along
$D_1, \ldots, D_n$, respectively, that leaves $U$ unknotted. Let $U'$ denote the result of $U$ after performing the crossings changes along $D_1, \ldots, D_{n-1}$ only, 
leaving $D_n$ intact. By (b) and (d), $U'\cup K_1 \cup\ldots\cup  K_{n-1}$ is the trivial link. Then, a crossing change of order $r_n>0$ 
along $D_n$ leaves $U'$ unknotted. By the argument in the proof of Theorem 2.2 in  \cite{kn:ak}, we conclude that  $L_n$ is the trivial link. This is a contradiction since $L_n$
is hyperbolic. \qed

\subsection{Proof of Theorem \ref{theo:mainknot}} 
Fix $n\in \N$ and let $L_n$ be a link as in Lemma \ref{lem:brunnian}.
We will consider the component $U$ of $L_n$ as a 0-framed unknot in $S^3$.
Given an
$n$-tuple of integers ${\bf r}:= (r_1, \ldots, r_n)$, with $r_i>2$,
let $M_n({\bf r})$ denote the 3-manifold obtained 
from $M_n$ as follows: For $1\leq i\leq n$, perform
Dehn filling with slope ${\displaystyle {{ 1 \over r_i} }}$ along $\partial {\eta(K_i)}$.
Let $K:={U}({\bf r})$ denote the image of $U$ in $M_n({\bf r})$.
Clearly, $M_n({\bf r}):=\overline{S^3\setminus \eta({U}({\bf r}))}$. Note, that since the linking number of $K_i$ and $U$ is zero
the framing on $K$ induced by that of $U$ is zero.
Thus $K$ is a 0-framed
knot in $S^3$ that is obtained from $U$ by  a generalized crossing of order
$r_1,\ldots, r_n$ along
$D_1, \ldots D_n$, respectively. By Lemma \ref{lem:brunnian} (e), $K$ is non-trivial.
\vskip 0.04in

{\it Claim:}
We have $U\equiv K({\rm mod}(r_i, \ 2))$, for every $1\leq i\leq n$.
\smallskip

{\it Proof of Claim:} Fix $1\leq i\leq n$. Consider $U(i)$ the knot obtained from
$U$  as follows: For every $1\leq j\neq i\leq n$, perform a generalized crossing change of order $r_j$
along $D_j$. By Lemma \ref{lem:brunnian} (d),  $U(i)$ is a (0-framed) trivial knot.
Since by construction, $K$ is obtained from $U:=U(i)$ by a generalized crossing change of order $r_i$
along $D_i$, we have $U\equiv K({\rm mod}(r_i, \ 2))$. This proves the claim. 
\smallskip  

Now by Corollary \ref{corol:colored}
we have that, if both $U$ and $K$ are given the 0-framing, then  $J_N(K, e_{r_i})=J_N(U, e_{r_i})$ , for $i=1,\ldots, n$. 
Suppose now that $U$ and $K$ are given any  framing $f \in \Z$. By Lemma 3.27 of \cite{kn:kim}, under the frame change from 0 to $f$
both of $J_N(K, t)$ and $J_N(U, t)$ are changed by the factor $t^{f(N^2-1)}$. Thus, the equation
$J_N(K, e_{r_i})=J_N(U, e_{r_i})$ remains true for every framing.
This proves the first part of the theorem.

Now we prove the remaining claims made in the statement of the theorem: By
Thurston's hyperbolic Dehn filling theorem (\cite{kn:thurston}), there is a constant 
$C_n:=C(L_n)>0$,  such that
if $r_1, \ldots, r_n> C_n$
then $M_n({\bf r})$ admits a complete hyperbolic structure
of finite volume. Thus $K:={U}({\bf r})$ is a hyperbolic knot.
By the proof of Thurston's theorem, the hyperbolic metric on
$M_n({\bf r})$ can be chosen so that it is arbitrarily close
to the metric of $M_n$, provided that the numbers
$r_i>>0$ are all sufficiently large.  Thus by choosing
the $r_i$'s large
we may ensure that the volume of $M_n({\bf r})$
is arbitrarily close to that of $M_n$. 
Since $\partial M_n$ has $n+1$ components, the interior of $M_n$
has $n+1$ cusps. By \cite{kn:adams}, we have ${\rm vol}(M_n)\geq ( n+1) v_3$,
where $v_3 (\approx 1.01494)$ is the volume of regular hyperbolic ideal tetrahedron.
Thus 
for $r_1, \ldots, r_n >>C_n$  we have

$${\rm vol}(S^3\setminus K)={\rm vol}(M_n)>n v_3> n. \eqno(2)$$

Now we turn our attention to closed 3-manifolds obtained by surgery of $S^3$ along 
a knot $K$ as above: Suppose that $M$ is a hyperbolic 3-manifold
obtained by ${\displaystyle {p\over q}}$-surgery on $K$. In Theorem 3.4 of \cite{dej} it is shown  that if $|q|>12$, then
${\rm vol}(M) \: \geq \: (1-{{127}\over {q^2}})^{3\over 2}
 {\rm vol}(S^3\setminus K)$. Combining this with (2) above we immediately
 obtain (1). This finishes the proof of the theorem. \qed


\subsection{Proof of Theorem \ref{theo:main}} Fix $n>0$ and $r_1, \ldots, r_n >2$
and let $K$ be
a knot as in Theorem \ref{theo:mainknot}. For a positive integer $q$,
let $M:=M_q(K)$ denote the 3-manifold obtained by surgery of
$S^3$ along $K$ with surgery slope $1\over q$. Let $L:= K^q$ denote the $q$-th cable of $K$ formed with 
the 0-framing as before. Similarly let $U^q$ denote the $q$-th cable of the unknot  $U$.
By Corollary \ref{corol:colored}, if both $L$ and $U^q$ are considered with $0$-framing,  we have 
$J_N(L, e_{r_i})=J_N(U^q, e_{r_i})$, for all $N\in \N$.
Now by Kirby calculus
$M$ can be obtained by surgery on the link $L:=K^q$ and where the framing on each of the
$q$ components is $1$ (see, for example, Figure 8 of \cite{gar}
for details.) 
By formula (1.9) on page 479 of \cite{kn:kim}, $\tau_{r_i}(M)$ is a linear combination of the values
$\{ J_N(L, e_{r_i}) \ | \  N< r_i \}$ with the linear coefficients depending
only on $r_i$ and the linking matrix of $L$. 
From this and our earlier observations,
the value of $\tau_{r_i}(M)$ remains the same if we replace $L$ with
$U^q$ and keep the same framings. But then the 3-manifold obtained by surgery on this later framed
link, which is the same as this obtained by $1\over q$-surgery on the unknot $U$, is clearly $S^3$.
Since, by \cite{kn:kim}, $\tau_{r}(S^3)=1$, for every $r>0$, we have
$$\tau_{r_i}(M)=\tau_{r_i}(S^3)=1, \ \ {\rm for}\ \ i=1, \ldots, m.$$ 
Next suppose that we have chosen  the values $r_1, \ldots, r_n >2$ large enough so that $K$ is hyperbolic.
By
Thurston's hyperbolic Dehn filling theorem, if $q>> 0$
then the 3-manifold $M$ is also hyperbolic.
We may, without loss of generality, assume that $|q|>12$. Now Theorem \ref{theo:mainknot}
implies that
${\rm vol}(M) \: \geq \: (1-{{127}\over {q^2}})^{3\over 2} n$, and for $q>>12$, we can assure that
 $${\rm vol}(M) \: \geq \: (1-{{127}\over {q^2}})^{3\over 2} n>{n \over 2}.$$ 
Next we show that the set
$A_K:=\{ M_q(K) \ \ | \ \  q\in\Z \}$,
contains infinitely many non-homeomorphic 3-manifolds. By \cite{kn:thurston}, for $q>>0$, we have
${\rm vol}( M_q(K))< {\rm vol}(S^3\setminus K)$ and ${\displaystyle{\lim_{q\to \infty}{\rm vol}( M_q(K))={\rm vol}(S^3\setminus K)}}$. Thus we can find a sequence  $\{M_k\}_{k\in \N}$
as claimed in the statement of Theorem \ref{theo:main}. \qed


\subsection { Proof of Corollary 1.3} 
We will use the notation and the setting established in the proofs Theorems \ref{theo:mainknot} and \ref{theo:main}:
Fix $n>1$ and choose $r'_1, \ldots, r'_n>>0$ as in the proof of Theorem \ref{theo:mainknot}, so that
for
$K:={U}({\bf r})$ we have ${\rm vol }(S^3\setminus K)>n$,
for every 
${\bf r}:= (r_1, \ldots, r_n)$ with $r_j \geq r'_j$.
For $i\in \N$ set
$${\bf r}_i:= (i+r_1, \ldots, r_n) \ \ {\rm and} \ \ K_i:={U}({\bf r}_i).$$
By the arguments in the proof of Theorem \ref{theo:main}, we have
$${\rm vol }(S^3\setminus K_i)>n \ \ {\rm and} \ \ K\equiv U({\rm mod}(m_i, \ 2)),$$
\no where $m_i:=i+r_1$. Thus Corollary \ref{corol:colored} and the argument in the proof of \ref{theo:main},
imply that there is a 3-manifold $M_i$ obtained by surgery along $K_i$ such that:
i) $\tau_{m_i}(M_i)=1$; and ii) $M_i$ is hyperbolic with ${\displaystyle {{\rm vol}(M_i)>{n \over 2}}}$.
\qed

\section{Concluding remarks}

1.\  Note that in the case that $N=r$ Corollary \ref{corol:colored}
simply states the well known fact $J_N(K, e_N)=J_N(J, e_N)=0$, for every knot $K$. Thus, in particular,
Corollary \ref{corol:colored} doesn't say anything  about the values of the colored Jones
polynomial that are relevant to the Volume Conjecture (\cite{kn:mumu}).
\smallskip

2. \ Relations between the volume and the $SU(2)$ invariants of 3-manifolds were
also studied by Kawauchi 
in \cite{kn:ka}. The main result in \cite{kn:ka}
implies the following: Given integers $R>0$ and $N>0$ there are
$2^N$ distinct, closed hyperbolic 3-manifolds that
share the same invariants $\tau_r$, for all levels $r< R$.
All of these $2^N$ manifolds, though, have the same volume
(in fact the same Chern-Simons invariants as well).  On the
other hand, Theorem \ref{theo:main} of this paper asserts the existence
of an infinite sequence of hyperbolic,  closed 3-manifolds, whose volumes form a strictly increasing
sequence, and all of which have the same $\tau_r$ for a finite set of levels.
\smallskip
\smallskip

{\bf Acknowledgment.} I thank Oliver Dasbach and Tim Cochran for helpful discussions and comments.
This research is supported in part by NSF grant DMS--0306995  and NSF--FRG grant DMS-0456155.

\end{document}